\documentclass[11pt]{amsart}

\usepackage{amsmath}
\usepackage{amssymb}

\newtheorem{theorem}{Theorem}[section]

\newtheorem{proposition}[theorem]{Proposition}
\newtheorem{corollary}[theorem]{Corollary}

\theoremstyle{definition}
\newtheorem{definition}[theorem]{Definition}

\newtheorem{question}[theorem]{Question}

\theoremstyle{remark}
\newtheorem{remark}[theorem]{Remark}

\newcount\skewfactor
\def\mathunderaccent#1#2 {\let\theaccent#1\skewfactor#2
\mathpalette\putaccentunder}
\def\putaccentunder#1#2{\oalign{$#1#2$\crcr\hidewidth
\vbox to.2ex{\hbox{$#1\skew\skewfactor\theaccent{}$}\vss}\hidewidth}}
\def\name{\mathunderaccent\tilde-3 }

\def\smallbox#1{\leavevmode\thinspace\hbox{\vrule\vtop{\vbox
   {\hrule\kern1pt\hbox{\vphantom{\tt/}\thinspace{\tt#1}\thinspace}}
   \kern1pt\hrule}\vrule}\thinspace}

\newcommand{\cf}{{\rm cf}}

\def\qedref#1{$\qed_{\reforiginal{#1}}$}

\title{Weak diamond and Galvin's property}
\author{Shimon Garti}
\address{Institute of Mathematics,
 The Hebrew University of Jerusalem,
 Jerusalem 91904, Israel}
\email{shimon.garty@mail.huji.ac.il}
\thanks{}


\begin{document}
\let\labeloriginal\label
\let\reforiginal\ref

\begin{abstract}
%%% Put abstract here:
Let $\kappa$ be an infinite cardinal, and $2^\kappa<\lambda\leq 2^{\kappa^+}$.
We prove that if there is a weak diamond on $\kappa^+$ then every $\{C_\alpha:\alpha<\lambda\}\subseteq\mathcal{D}_{\kappa^+}$ satisfies Galvin's property. On the other hand, Galvin's property is consistent with the failure of the weak diamond (and even with Martin's axiom in the case of $\aleph_1$). We derive some consequences about weakly inaccessible cardinals. We also prove that the negation of a similar property follows from the proper forcing axiom.
\end{abstract}

\maketitle

% start document here:
\newpage

\section{Introduction}

Assume $\kappa$ is a regular cardinal, and let $\mathcal{D}_\kappa$ be the club filter on $\kappa$.
Being $\kappa$-complete, the intersection of $\theta$-many members of $\mathcal{D}_\kappa$ contains a club subset of $\kappa$, whenever $\theta<\kappa$. However, the limitation $\theta<\kappa$ is essential. Indeed, every end-segment of $\kappa$ is a club, and the intersection of all the end-segments is empty. On the other hand, it is easy to define a collection of $2^\kappa$ different clubs in $\kappa$ such that the intersection of which contains a club (e.g., take a club $C$ such that $\kappa\setminus C$ is of size $\kappa$, and define all the clubs as $C\cup A$ for some $A\subseteq\kappa\setminus C$).

The interesting question is, therefore, whether for any collection of $\lambda$ clubs one can find a sub-collection of cardinality $\kappa$ such that the intersection of which is a club. We shall work with the club filter $\mathcal{D}_\kappa$. An element of $\mathcal{D}_\kappa$ need not be closed, but it contains a club.
Notice that in the specific case of the club filter, if the intersection contains an unbounded set of $\kappa$ then it contains a club (just take the closure of the unbounded set in the order topology).
Galvin's property is exactly this situation. Galvin proved that if $\kappa=\kappa^{< \kappa}$ then any family $\{C_\alpha:\alpha<\kappa^+\}\subseteq\mathcal{D}_{\kappa}$ contains a subfamily of cardinality $\kappa$ with unbounded intersection (and hence also a club included in the intersection).

The proof appears in \cite{MR0369081}, and attributed to Fred Galvin. On one hand, the proof is phrased for the specific case of $\kappa=\aleph_1$, and then the pertinent assumption is $2^{\aleph_0}=\aleph_1$. On the other hand, the proof there is more general as it applies to every normal filter (rather than the club filter). A natural question is whether the assumption $\kappa=\kappa^{< \kappa}$ is required.
Abraham and Shelah proved in \cite{MR830084} that some assumption must be made, as the negation of Galvin's property is consistent.

We shall see that one can use a weaker assumption in order to get a general form of Galvin's property. The conclusion would be that in the model of Abraham-Shelah we have $2^\kappa=2^{\kappa^+}$. It means that the principle of weak diamond (from \cite{MR0469756}) fails. In other words, if the weak diamond holds then for some $\lambda\leq 2^{\kappa^+}$ we can prove Galvin's property. As in the original proof, the result holds true for every normal filter on $\kappa$.

We phrase the definition of the weak diamond:

\begin{definition}
\label{wweakdiamond} The Devlin-Shelah weak diamond. \newline
Let $\kappa$ be an infinite cardinal. \newline
The weak diamond on $\kappa^+$ (denoted by $\Phi_{\kappa^+}$) is the following principle: \newline
For every function $F:{}^{<\kappa^+} 2\rightarrow 2$ there exists a function $g\in{}^{\kappa^+} 2$ such that $\{\alpha\in\kappa^+: F(f\upharpoonright \alpha)=g(\alpha)\}$ is a stationary subset of $\kappa^+$ whenever $f\in{}^{\kappa^+} 2$.
\end{definition}

The main content of \cite{MR0469756} is the proof that the weak diamond (on $\aleph_1$) follows by the assumption that $2^{\aleph_0}<2^{\aleph_1}$. However, the opposite direction is also true, i.e. if there is a weak diamond on $\aleph_1$ then $2^{\aleph_0}<2^{\aleph_1}$. The proof is attributed to Uri Abraham, in \cite{MR1623206}. Both directions generalize, verbatim, to every successor cardinal $\kappa^+$ in lieu of $\aleph_1$. We spell out the easy direction (see \cite{MR1623206} pp. 943-944 for $\aleph_1$):

\begin{proposition}
\label{oooo} If $\Phi_{\kappa^+}$ then $2^\kappa<2^{\kappa^+}$.
\end{proposition}

\par\noindent\emph{Proof}. \newline
Assume towards contradiction that $\Phi_{\kappa^+}$ holds but $2^\kappa=2^{\kappa^+}$. Choose a bijection $b:2^\kappa\rightarrow 2^{\kappa^+}$. We define $F:{}^{<\kappa^+}2\rightarrow 2$ as follows. If $\eta\in {}^\alpha 2$ for some $\alpha\in\kappa$ then $F(\eta)=0$ and if $\eta\in {}^\alpha 2$ for some $\alpha\in[\kappa,\kappa^+)$ then $F(\eta)=[b(\eta\upharpoonright\kappa)](\alpha)$.

We claim that $F$ cannot be predicted by any $g\in {}^{\kappa^+}2$, so assume to the contrary that some $g\in {}^{\kappa^+}2$ exemplifies the weak diamond for $F$. For every $\alpha\in[\kappa,\kappa^+)$ set $h(\alpha)=1-g(\alpha)$ and for every $\alpha<\kappa$ define $h(\alpha)=0$. Since $h\in {}^{\kappa^+}2$, there is some $t\in {}^\kappa 2$ such that $b(t)=h$. Extend $t$ to a mapping $f\in {}^{\kappa^+}2$ by letting $f(\alpha)=0$ for every $\alpha\in[\kappa,\kappa^+)$ and $f(\alpha)=t(\alpha)$ for $\alpha<\kappa$.

Now for every $\alpha\in[\kappa,\kappa^+)$ we have $F(f\upharpoonright\alpha)= [b(f\upharpoonright\kappa)](\alpha)= b(t)(\alpha)=h(\alpha)\neq g(\alpha)$. It follows that $F(f\upharpoonright\alpha)=g(\alpha)$ is possible only for $\alpha<\kappa$ which is not a stationary set in $\kappa^+$, so $\Phi_{\kappa^+}$ fails, a contradiction.

\hfill \qedref{oooo}

The above discussion focuses on successor cardinals. Similar results about Galvin's property can be proved for inaccessible cardinals. The interesting case is a weakly inaccessible $\lambda$ so that $\lambda^{<\lambda}>\lambda$. We shall see that the weak diamond on an inaccessible $\lambda$ implies Galvin's property for subsets of $\lambda$. This gives a partial answer to Question 1.28 from \cite{MR2210145}.

``Consistency... I advocate consistency in all things" (\cite{ja}, p. 33).
The last topic is Galvin's property and forcing axioms. It is shown here that Martin's axiom is consistent with Galvin's property, although it is consistent also with its negation. We suspect that stronger axioms (e.g. Martin's maximum or the PFA) imply the failure of Galvin's property. Here we prove a partial result in this direction under the PFA.

We conclude this section with some basic facts and definitions. We use the Jerusalem notation in our forcing, so $p\leq q$ means that $q$ extends $p$, and hence a generic set $G$ is \emph{downward} closed. The properness of a forcing notion $\mathbb{Q}$ can be described by the existence of $(M,\mathbb{Q})$-generic conditions for every suitable elementary submodel $M$.
The \emph{Proper forcing axiom} says that for every proper forcing $\mathbb{P}$, if $\mathcal{D}$ is a collection of dense subsets of $\mathbb{P}$ and $|\mathcal{D}|\leq\aleph_1$ then there is a generic subset $G\subseteq\mathbb{P}$ so that $G\cap D\neq\emptyset$ for every $D\in\mathcal{D}$.

The \emph{Bounded proper forcing axiom} is a relative of the PFA, first defined in \cite{MR1324501}. It is similar to the PFA, but promises a generic set only for a proper forcing notion $\mathbb{P}$ with an extra requirement that every maximal antichain $\mathcal{A}$ of $\mathbb{P}$ is of size $\aleph_1$. Our forcing, below, satisfies this requirement.

We shall prove the failure of a weak form of Galvin's property under the PFA. For this, we will use the forcing notion of Baumgartner, \cite{MR776640}, which adds a club subset to $\omega_1$ through finite conditions. The pertinent theorem is the following:

\begin{theorem}
\label{ppffaa} (PFA). \newline
Assume $\langle t_\beta:\beta<\omega_1\rangle$ is a sequence of infinite subsets of $\omega_1$. Then there exists a closed unbounded set $C\subseteq\aleph_1$ so that for every $\beta<\omega_1$, $\neg(t_\beta\subseteq C)$.
\end{theorem}

The original proof of the theorem is ensconced in \cite{MR776640}, Theorem 3.4. Stronger results in this vein can be found in \cite{MR1623206}.
We refer to \cite{MR2768684} for a comprehensive account of proper forcing, as well as \cite{MR1623206}.

I wish to thank the referee for the corrections and improvements. I also thank Uri Abraham for a very helpful discussion concerning the contents of this paper, and Martin Goldstern for several clever remarks.

\newpage

\section{Weak diamond, MA and PFA}

We commence with a generalization of Galvin's proof:

\begin{theorem}
\label{mt} Large club intersection. \newline
Assume $2^\kappa<2^{\kappa^+}$, and $\lambda\in(2^\kappa,2^{\kappa^+}]$. \newline
For every family $\{C_\alpha:\alpha<\lambda\}\subseteq\mathcal{D}_{\kappa^+}$ there is a subfamily of cardinality $\kappa^+$ whose intersection belongs to $\mathcal{D}_{\kappa^+}$. \newline
The same holds for an inaccessible cardinal $\kappa$, when $\kappa=\kappa^{<\kappa}$ and $\lambda\in[\kappa^+,2^\kappa]$. It means that for every family $\{C_\alpha:\alpha<\lambda\}\subseteq\mathcal{D}_\kappa$ one can find a subfamily of size $\kappa$ whose intersection is a member of $\mathcal{D}_\kappa$.
\end{theorem}

\par\noindent\emph{Proof}. \newline
Given a family of more than $2^\kappa$ members of $\mathcal{D}_{\kappa^+}$, we can concentrate on $(2^\kappa)^+$-many of them, so without loss of generality, $\lambda=(2^\kappa)^+$ (and hence a regular cardinal). For every $\alpha<\lambda$ and $\varepsilon<\kappa^+$, set:

$$
H_{\alpha\varepsilon}=\{\beta<\lambda: C_\alpha\cap\varepsilon = C_\beta\cap\varepsilon\}.
$$

The basic observation of Galvin (generalized here) is that for some $\alpha<\lambda$ we have $|H_{\alpha\varepsilon}|=\lambda$ for every $\varepsilon<\kappa^+$. For suppose the contrary, and for every $\alpha<\lambda$ choose an ordinal $\varepsilon_\alpha<\kappa^+$ so that $|H_{\alpha\varepsilon_\alpha}|<\lambda$. Since $\lambda=\cf(\lambda)>2^\kappa\geq\kappa^+$, there is a subset $S\in[\lambda]^\lambda$ and an ordinal $\varepsilon<\kappa^+$ such that $\alpha\in S\Rightarrow \varepsilon_\alpha\equiv\varepsilon$.

As $\varepsilon<\kappa^+$ we know that $|\mathcal{P}(\varepsilon)|\leq 2^\kappa<\lambda$. Observe that the set $H_{\alpha\varepsilon}$ is determined absolutely by $C_\alpha\cap\varepsilon\subseteq\varepsilon$, so we have among the family $\{H_{\alpha\varepsilon}:\alpha\in S\}$ less than $\lambda$-many different sets. By the assumption towards contradiction, the size of each member is less than $\lambda$, and hence $|\bigcup\limits_{\alpha\in S}H_{\alpha\varepsilon}|<\lambda$. On the other hand, for every $\alpha\in S$ we have trivially $\alpha\in H_{\alpha\varepsilon}$, and since $|S|=\lambda$ we conclude that $|\bigcup\limits_{\alpha\in S}H_{\alpha\varepsilon}|=\lambda$ as well, a contradiction.

Fix an ordinal $\alpha<\lambda$ so that $|H_{\alpha\varepsilon}|=\lambda$ for every $\varepsilon<\kappa^+$. By induction on $\varepsilon<\kappa^+$ we choose an ordinal $\zeta_\varepsilon$ such that $\zeta_\varepsilon\in H_{\alpha\varepsilon+1}\setminus \{\zeta_\eta:\eta<\varepsilon\}$. There is no problem in carrying the induction, since $|H_{\alpha\varepsilon+1}|=\lambda$ at every stage.

We claim that the collection $\{C_{\zeta_\varepsilon}:\varepsilon<\kappa^+\}$ exemplifies Galvin's property. For this, define $C=\bigcap\limits_{\varepsilon<\kappa^+} C_{\zeta_\varepsilon}$, and we shall prove that $C$ is a club subset of $\kappa^+$. We intend to show that $C\cup (\kappa^+\setminus C_\alpha)\in\mathcal{D}_{\kappa^+}$, and since $\kappa^+\setminus C_\alpha$ belongs to the non-stationary ideal of $\kappa^+$ we are getting $C\in \mathcal{D}_{\kappa^+}$. In order to prove this assertion, we shall prove that:

$$
C\cup (\kappa^+\setminus C_\alpha)\supseteq \Delta\{C_{\zeta_\varepsilon}:\varepsilon<\kappa^+\}\in\mathcal{D}_{\kappa^+}.
$$

Assume that $\beta\in \Delta\{C_{\zeta_\varepsilon}:\varepsilon<\kappa^+\}$. If $\beta\in \kappa^+\setminus C_\alpha$ we are done, so assume $\beta\in C_\alpha$. We have to show that $\beta\in C_{\zeta_\varepsilon}$ for every $\varepsilon<\kappa^+$. Indeed, for $\varepsilon<\beta$ the membership $\beta\in C_{\zeta_\varepsilon}$ follows from the definition of the diagonal intersection. If $\varepsilon\geq\beta$ then $\beta<\varepsilon+1$ and hence $\beta\in C_\alpha\cap(\varepsilon+1)=C_{\zeta_\varepsilon}\cap(\varepsilon+1)$ since $\zeta_\varepsilon\in H_{\alpha\varepsilon+1}$, so in particular $\beta\in C_{\zeta_\varepsilon}$. It means that $\beta\in \bigcap\limits_{\varepsilon<\kappa^+} C_{\zeta_\varepsilon}=C$.
The proof in the inaccessible case under the assumption $\kappa=\kappa^{<\kappa}$ is just the same (upon replacing $\kappa^+$ by $\kappa$), so we are done.

\hfill \qedref{mt}

If $\lambda>2^{\kappa^+}$ then it is virtually trivial that every collection of $\lambda$ clubs of $\kappa^+$ contains a subfamily of size $\kappa^+$ whose intersection is a club. If $\kappa=\theta^+$ and $\theta>\aleph_0$ then the assumption that $\kappa=\kappa^{<\kappa}$ means that $2^\theta=\theta^+$, in which case $\Diamond_{\theta^+}$ holds.
The interesting consequence of the above theorem is that Galvin's property for cardinals below $2^{\kappa^+}$ is strongly tied to the weak diamond principle:

\begin{corollary}
\label{mc} Weak diamond and large club intersection. \newline
Assume $\Phi_{\kappa^+}$. \newline
There exists some $\lambda\leq 2^{\kappa^+}$ such that every collection of $\lambda$ clubs of $\kappa^+$ contains a subfamily of size $\kappa^+$ whose intersection contains a club.
\end{corollary}

\hfill \qedref{mc}

It is worth noticing that the above claim and the original theorem of Galvin can be incorporated under the same canopy in the following manner:
If $\kappa=\theta^+$ and $\{C_\alpha:\alpha<(2^\theta)^+\}\subseteq\mathcal{D}_\kappa$, then there exists a subfamily of $\kappa$-many sets whose intersection belongs to $\mathcal{D}_\kappa$.
Indeed, if $2^\theta=2^\kappa$ then the proposition follows from the pigeon-hole principle. If $2^\theta<2^\kappa$ then it follows from Theorem \ref{mt}. If, in addition, $2^\theta=\kappa$, then this is the theorem of Galvin.
So as a matter of fact, Galvin's property is a general theorem, which is not confined to the local instance of the continuum hypothesis. It holds in ZFC, provided that the pertinent $\lambda$ (i.e., the number of club sets) is designated properly. We collect some useful comments concerning the proof:

\begin{remark} Trifles.
\begin{enumerate}
\item [$(a)$] In order to choose a new ordinal $\zeta_\varepsilon$ one needs merely that $|H_{\alpha\varepsilon+1}|\geq\kappa^+$ (rather than $|H_{\alpha\varepsilon+1}|=\lambda$).
\item [$(b)$] The equality $C_\beta\cap\varepsilon=C_\alpha\cap\varepsilon$ can be replaced by $C_\beta\cap\varepsilon\supseteq C_\alpha\cap\varepsilon$.
\item [$(c)$] The proviso $|H_{\alpha\varepsilon}|=\lambda$ \emph{for every} $\varepsilon<\kappa^+$ can be replaced by an unbounded set of $\varepsilon$-s below $\kappa^+$.
\item [$(d)$] The filter $\mathcal{D}_{\kappa^+}$ can be replaced by any normal filter on $\kappa^+$.
\item [$(e)$] If $\kappa$ is measurable then $\binom{\kappa^+}{\kappa}\rightarrow\binom{\kappa}{\kappa}^{1,1}_2$ and much more. If $D$ is a normal ultrafilter on $\kappa$ then the last part of the above theorem shows that the righthand component of the monochromatic product can be a member of $D$, so we may write $\binom{\kappa^+}{\kappa}\rightarrow_D \binom{\kappa}{\kappa}^{1,1}_2$. Indeed, the member of $D$ which comes from Galvin's property can serve as this component.
\end{enumerate}
\end{remark}

Concerning $(b)$, it is possible to render the proof of Theorem \ref{mt} with a diamond sequence (instead of $2^\kappa=\kappa^+$). Part $(b)$ of the remark suggests a weaker possibility:

\begin{question}
\label{q0}
Does $\clubsuit_{\kappa^+}$ imply Galvin's property on $\kappa^+$?
\end{question}

A natural question arises: Is it sufficient to violate the weak diamond in order to eliminate Galvin's property? The following theorem shows that Galvin's property is consistent with the negation of the weak diamond. Moreover, Martin's axiom is a strong environment of the failure of weak diamonds (as $MA+2^{\aleph_0}=\lambda\Rightarrow 2^{\aleph_0}=2^\theta$ for every $\theta<\lambda$), but we shall see that Galvin's property is independent over Martin's axiom.

\begin{theorem}
\label{nnegg} Galvin's property and $\neg\Phi$. \newline
Let $\kappa$ be an infinite cardinal. \newline
It is consistent that for some $\lambda\leq 2^{\kappa^+}$ every collection $\{C_\alpha:\alpha<\lambda\}\subseteq\mathcal{D}_{\kappa^+}$ contains a subfamily of $\kappa^+$-many members whose intersection belongs to $\mathcal{D}_{\kappa^+}$ while $\neg\Phi_{\kappa^+}$. \newline
Moreover, the above Galvin's property is consistent even with Martin's axiom, and even for $\lambda=\kappa^{++}$.
\end{theorem}

\par\noindent\emph{Proof}. \newline
We begin with $2^\kappa<2^{\kappa^+}$ in the ground model (and assuming further that $2^\kappa=\kappa^+$ we can choose $\lambda$ below to be $\kappa^{++}$). Let $\mathbb{P}$ be any $\kappa^+$-cc forcing notion such that if $G\subseteq\mathbb{P}$ is generic then $V[G]\models 2^\kappa=2^{\kappa^+}$. Since $\mathbb{P}$ is $\kappa^+$-cc we claim that $\mathcal{D}_{\kappa^+}^{V[G]}$ is generated by $\mathcal{D}_{\kappa^+}^V$, i.e. for every $C\in\mathcal{D}_{\kappa^+}^{V[G]}$ there exists $C'\in \mathcal{D}_{\kappa^+}^V$ so that $C'\subseteq C$.

For proving this fact, let $f\in V[G]$ be the enumeration of the members of $C$ in increasing order, so $f$ is a normal function. We can choose a condition $p\in G$ such that $p \Vdash \name{f}:\check{\kappa}\rightarrow \name{C}, \name{f}$ is normal. We define, by induction on $\kappa^+$, a function $g:\kappa^+\rightarrow\kappa^+$ as follows:
if $\alpha=0$ then $g(\alpha)=0$, and if $\alpha$ is a limit ordinal then $g(\alpha)=\bigcup\limits_{\beta<\alpha} g(\beta)$. If $\alpha=\beta+1$ we choose $\gamma,\delta<\kappa^+$ such that $p \Vdash g(\beta)<\name{f}(\check{\gamma})<\check{\delta}$. Here we use the chain condition, as the value of $\name{f}(\check{\gamma})$ is determined by an antichain of size less than $\kappa^+$, and hence we can find an ordinal $\delta$ above all possible values and below $\kappa^+$. We define $g(\alpha)=\delta$. Notice that $g\in V$.

Let $\eta$ be any limit ordinal below $\kappa^+$.
Since $p$ forces that $\name{f}$ is a normal function, and by the interlaced nature of $f$ and $g$, we know that $p \Vdash g(\eta)\in \name{C}$. So we can define (in $V$) the set $C'=\{g(\eta):\eta$ is a limit ordinal$\}$. It follows that $C'\in V$ as claimed.

Choose any $\lambda\in[(2^\kappa)^+,2^{\kappa^+}]$. Given any family $\{C_\alpha:\alpha<\lambda\}\subseteq\mathcal{D}_{\kappa^+}^{V[G]}$ we choose for every $\alpha<\lambda$ some $C'_\alpha\in \mathcal{D}_{\kappa^+}^V$ such that $C'_\alpha\subseteq C_\alpha$.
Using once again the $\kappa^+$-cc, we may assume without loss of generality that $\{C'_\alpha:\alpha<\lambda\}\subseteq\mathcal{D}_{\kappa^+}^V$. More precisely, we have to assume that $\lambda$ is regular (for this we may assume that $\lambda=(2^\kappa)^+$) and then we can argue that some unbounded family of $\{C'_\alpha:\alpha<\lambda\}$ belongs to $\mathcal{D}_{\kappa^+}^V$.
The argument can be justified by the property that we proved above, since we can find in ${\rm V}$ an unbounded set of indices $\alpha$ below $\lambda$ for which $C'_\alpha\subseteq C_\alpha$ and $C'_\alpha\in \mathcal{D}_{\kappa^+}^V$. As $\lambda=\cf(\lambda)$ we can rename and get the desired collection $\{C'_\alpha:\alpha<\lambda\}$.
However, the ground model satisfies Galvin's property with respect to $\lambda$, hence there is a subfamily $\{C'_{\alpha_\varepsilon}:\varepsilon<\kappa^+\}$ for which $C'=\bigcap\limits_{\varepsilon<\kappa^+}C'_{\alpha_\varepsilon}\in \mathcal{D}_{\kappa^+}^V$.
Notice that $C'\subseteq \bigcap\limits_{\varepsilon<\kappa^+}C_{\alpha_\varepsilon}$ as well, and since $C'\in\mathcal{D}_{\kappa^+}^V \subseteq \mathcal{D}_{\kappa^+}^{V[G]}$ we are done.

For the additional part of the claim, recall that if $\mu=\mu^{<\mu}>\aleph_0$ then there is an $\aleph_1$-cc forcing notion $\mathbb{P}$ such that if $G\subseteq\mathbb{P}$ is generic then Martin's axiom holds in $V[G]$ (see, e.g. \cite{MR1623206}).
Begin with a ground model in which $2^{\aleph_0}=\aleph_1$, e.g. $V=L$. Notice that Galvin's property holds for $\aleph_1$ in the ground model. Force with $\mathbb{P}$ in order to get Martin's axiom and $2^{\aleph_0}=\mu$ in $V[G]$. Since $\mathbb{P}$ is $\aleph_1$-cc it preserves $\mathcal{D}_{\aleph_1}$ and hence Galvin's property for $\aleph_1$ holds in $V[G]$ as well.

\hfill \qedref{nnegg}

\begin{remark} More trifles.
\begin{enumerate}
\item [$(a)$] The above theorem gives the consistency of Galvin's property for $\lambda=\kappa^{++}$ while $2^\kappa=2^{\kappa^+}$ is arbitrarily large.
\item [$(b)$] It clarifies the fact that the forcing notion in \cite{MR830084} is not $\kappa^+$-cc (and cannot be).
\item [$(c)$] The negation of Galvin's property is also consistent with Martin's axiom, see \cite{MR830084}.
\end{enumerate}
\end{remark}

Concerning the last remark, although Martin's axiom implies $\neg(\Phi_{\aleph_1})$ there is a difference between the negation of the weak diamond and Martin's axiom. It seems that Martin's axiom is a meaningful step towards the negation of Galvin's property, and we may ask what happens under stronger forcing assumptions:

\begin{question}
\label{q1}
Is it consistent that Galvin's property (for $\aleph_1$) holds under the PFA?
\end{question}

The general impression is that Galvin's property fails under the PFA. We shall prove somewhat weaker assertion.
We deal with clubs of $\aleph_1$, and we indicate that the pertinent number of clubs under the PFA is $\aleph_2$ as $2^{\aleph_1}=\aleph_2$.
In the theorem below we first prove a claim which follows from $2^{\aleph_1}=\aleph_2$, and then we prove a stronger claim which follows from the PFA. Notice that in the second part we get finite intersections (a stronger property than bounded intersections in $\aleph_1$):

\begin{theorem}
\label{ppffa1} PFA and club intersection. \begin{enumerate}
\item [$(a)$] Under the assumption $2^{\aleph_1}=\aleph_2$, there exists a collection $\{C_\alpha:\alpha<\omega_2\}\subseteq \mathcal{D}_{\aleph_1}$ such that the intersection of every $\aleph_2$ of them is bounded. A parallel assertion holds for every regular uncountable cardinal $\kappa$, under the assumption $2^\kappa=\kappa^+$.
\item [$(b)$] Under the PFA, there exists a collection $\{C_\alpha:\alpha<\omega_2\}\subseteq \mathcal{D}_{\aleph_1}$ such that the intersection of every $\aleph_2$ of them is finite. Moreover, the PFA can be replaced by the BPFA.
\end{enumerate}
\end{theorem}

\par\noindent\emph{Proof}. \newline
For part $(a)$ we first show that if $\{A_\beta:\beta<\aleph_1\}$ is any family of unbounded subsets of $\aleph_1$ then there is a club $C\subseteq\aleph_1$ for which $A_\beta\nsubseteq C$ for every $\beta<\aleph_1$. In order to prove this, we define an unbounded set $C'=\{\gamma_\beta:\beta<\aleph_1\} \subseteq\aleph_1$ by induction. For $\beta=0$ let $a_0$ be the first member of $A_0$ and we choose an ordinal $\gamma_0>a_0$. For $\beta>0$, if $\gamma_\alpha$ has been chosen for every $\alpha<\beta$ then we choose $a_\beta\in A_\beta$ such that $a_\beta>\bigcup\{\gamma_\alpha:\alpha<\beta\}$ (this can be done since $A_\beta$ is unbounded in $\aleph_1$) and then we choose an ordinal $\gamma_\beta>a_\beta$. Let $C'$ be $\{\gamma_\beta:\beta<\aleph_1\}$.
Our club $C$ is the closure of $C'$. It follows from the construction that $a_\beta\notin C$ for every $\beta\in\aleph_1$, and hence $A_\beta\nsubseteq C$ for every $\beta<\aleph_1$.

Now let $\{A_\beta:\beta<\omega_2\}$ enumerate all the members of $[\aleph_1]^{\aleph_1}$. By induction on $\alpha<\omega_2$ we choose a club $C_\alpha\subseteq\aleph_1$ so that $A_\beta\nsubseteq C_\alpha$ for every $\beta<\alpha$. This can be done by the above paragraph, since for every $\alpha<\omega_2$ there are at most $\aleph_1$-many sets in the collection $\{A_\beta:\beta<\alpha\}$. Finally, our collection is $\{C_\alpha:\alpha<\omega_2\}$ and we claim that it satisfies the required property. Indeed, if $\{C_{\alpha_\varepsilon}: \varepsilon<\omega_2\}\subseteq \{C_\alpha:\alpha<\omega_2\}$ and $A\subseteq \bigcap\limits_{\varepsilon<\omega_2} C_{\alpha_\varepsilon}$ is unbounded in $\aleph_1$, then $A=A_\beta$ for some $\beta<\omega_2$. But then for a large enough $\varepsilon<\omega_2$ we have $\beta<\alpha_\varepsilon$ and hence $A_\beta\nsubseteq C_{\alpha_\varepsilon}$, a contradiction.

We turn to the second assertion of the theorem. The proof is similar, but we have now a better running away property.
Fix an enumeration $\{t_\beta:\beta<\omega_2\}$ of all the members of $[\aleph_1]^{\aleph_0}$ (recall that $2^{\aleph_0}=2^{\aleph_1}=\aleph_2$ under the PFA, as proved in \cite{MR1174395}). We choose, by induction on $\alpha<\omega_2$, a club $C_\alpha$ in $\aleph_1$ such that:

$$
\forall\beta<\alpha, t_\beta\nsubseteq C_\alpha.
$$

This can be done by Theorem \ref{ppffaa}. Let $\{C_{\alpha_\varepsilon}: \varepsilon<\omega_2\}\subseteq \{C_\alpha:\alpha<\omega_2\}$. Assume towards contradiction that there is an infinite set $t\subseteq\omega_1$ so that $t\subseteq \bigcap\limits_{\varepsilon<\omega_2} C_{\alpha_\varepsilon}$. Without loss of generality $t$ is countable, so we can pick an ordinal $\beta<\omega_2$ such that $t\equiv t_\beta$.
Choose a large enough $\varepsilon<\omega_2$ such that $\beta<\alpha_\varepsilon$. By the definition of $C_{\alpha_\varepsilon}$ it follows that $t=t_\beta\nsubseteq C_{\alpha_\varepsilon}$, a contradiction, so we are done.

For the additional part of the statement let $\mathbb{Q}$ be the forcing notion of Theorem \ref{ppffaa}, and recall that $|\mathbb{Q}|=\aleph_1$. Hence, if $\mathcal{A}\subseteq\mathbb{Q}$ is a maximal antichain then $|\mathcal{A}|\leq\aleph_1$, and actually it is exactly $\aleph_1$.
In addition, $2^{\aleph_0}=\aleph_2$ under the BPFA, as proved in \cite{MR2151584}. Consequently, the result of Baumgartner in Theorem \ref{ppffaa} follows from the BPFA only.
It means that the BPFA yields the same conclusion, and the proof of the theorem is accomplished.

\hfill \qedref{ppffa1}

We turn to inaccessible cardinals.
Let $\lambda$ be a regular limit cardinal.
If $\lambda$ is strongly inaccessible then $\lambda=\lambda^{<\lambda}$, and this case is covered in Theorem \ref{mt}. Therefore, the interesting assumption is $\lambda<\lambda^{<\lambda}$, in which case $\lambda$ is weakly inaccessible.
We have seen, for successor cardinals, that the weak diamond implies Galvin's property. The following theorem shows that the same implication holds for limit cardinals as well.

\begin{theorem}
\label{inac} Galvin's property for weakly inaccessible cardinals. \newline
Let $\kappa$ be a weakly inaccessible cardinal.
\begin{enumerate}
\item [$(a)$] If $2^{<\kappa}<\mu\leq 2^\kappa$ and $\{C_\alpha:\alpha<\mu\}\subseteq\mathcal{D}_\kappa$ then there exists a subfamily of size $\kappa$ whose intersection belongs to $\mathcal{D}_\kappa$.
\item [$(b)$] If $\Phi_\kappa$ then $2^{<\kappa}<2^\kappa$ and hence Galvin's property holds with respect to every $\mu\in(2^{<\kappa},2^\kappa]$.
\item [$(c)$] It is consistent that $\Phi_\kappa$ fails while Galvin's property holds.
\end{enumerate}
\end{theorem}

\par\noindent\emph{Proof}. \newline
We commence with part $(a)$.
We may assume that $\mu=(2^{<\kappa})^+$, as we can choose a subfamily of the given sets. As usual, we define $H_{\alpha\varepsilon}=\{\beta<\mu: C_\alpha\cap\varepsilon = C_\beta\cap\varepsilon\}$ for every $\alpha<\mu$ and every $\varepsilon<\kappa$. By the same token, for some (and actually, many) $\alpha<\mu$ we know that $|H_{\alpha\varepsilon}|=\mu$ for every $\varepsilon<\kappa$. Here we use the fact that $|\mathcal{P}(\varepsilon)|\leq 2^{<\kappa}<\mu=\cf(\mu)$.

We choose an ordinal $\alpha<\mu$ with the above property. By induction (this time on $\varepsilon<\kappa$) we choose $\zeta_\varepsilon$ such that $\zeta_\varepsilon\in H_{\alpha\varepsilon+1}\setminus \{\zeta_\eta:\eta<\varepsilon\}$. We claim that $\{C_{\zeta_\varepsilon}:\varepsilon<\kappa\}$ is as required. Indeed, the intersection $C=\bigcap\limits_{\varepsilon<\kappa} C_{\zeta_\varepsilon}$ contains $\Delta\{C_{\zeta_\varepsilon}:\varepsilon<\kappa\}\cap C_\alpha$ and hence belongs to $\mathcal{D}_\kappa$.

For part $(b)$, assume towards contradiction that $\Phi_\kappa$ holds but $2^{<\kappa}=2^\kappa$. Fix a function $b:2^{<\kappa}\rightarrow 2^\kappa$ which is onto, and satisfies the following requirement: For every $\alpha<\kappa$ and every $t\in{}^\alpha 2$, if there is an ordinal $\delta<\alpha$ such that $\varepsilon\in[\delta,\alpha)\Rightarrow t(\varepsilon)=0$ then $b(t)=b(t\upharpoonright\delta)$.
We wish to define a coloring $F:{}^{<\kappa}2\rightarrow 2$, so assume $\alpha<\kappa$ and $\eta\in{}^\alpha 2$. Set $F(\eta)=[b(\eta)](\alpha)$.

By $\Phi_\kappa$ we can choose a function $g\in{}^\kappa 2$ which predicts $F$. Let $h\in{}^\kappa 2$ be the opposite function, i.e. $h(\alpha)=1-g(\alpha)$. Pick any function $t\in{}^{<\kappa}2$ for which $b(t)=h$. Let $\delta<\kappa$ be such that $t\in{}^\delta 2$. We define $f\in{}^\kappa 2$ as an extension of $t$ as follows. If $\alpha<\delta$ then $f(\alpha)=t(\alpha)$ and if $\alpha\geq\delta$ then $f(\alpha)=0$. Observe that $b(f\upharpoonright\alpha)=b(f\upharpoonright\delta)$ for every $\alpha\in(\delta,\kappa)$. Consequently, $F(f\upharpoonright\alpha)=[b(f\upharpoonright\alpha)](\alpha) =[b(f\upharpoonright\delta)](\alpha)=[b(t)](\alpha) =h(\alpha)\neq g(\alpha)$, a contradiction.

Part $(c)$ is proved as in the successor case. We begin with a weakly inaccessible $\kappa$ such that $2^{<\kappa}<2^\kappa$, and we force $2^\theta=2^\kappa$ for some $\theta<\kappa$ by a forcing notion which preserves $\mathcal{D}_\kappa$, and the weak inaccessibility of $\kappa$. Due to $(b)$, the weak diamond fails on $\kappa$ but Galvin's property holds since the new $\mathcal{D}_\kappa$ is generated by the old one.

\hfill \qedref{inac}

The weak diamond implies Galvin's property, so this property can be viewed as a feeble form of the weak diamond. Assume $\lambda$ is weakly inaccessible. If $\lambda=\lambda^{<\lambda}=2^\kappa$ for some $\kappa<\lambda$ then $\Phi_\lambda$ holds, as proved in \cite{MR1901038}. Suppose $\lambda<\lambda^{<\lambda}$, moreover $\{2^\theta:\theta<\lambda\}$ is not eventually constant. It means that unboundedly many instances of the weak diamond hold below $\lambda$. On this ground, Question 1.28 from \cite{MR2210145} is the following natural problem:
\begin{enumerate}
\item [$(\aleph)$] Does $\lambda$ carry the weak diamond?
\item [$(\beth)$] Can other consequences be proved on $\lambda$?
\end{enumerate}

We can answer partially, by the following:

\begin{corollary}
\label{cc1}
Assume $\lambda$ is weakly inaccessible and $\{2^\theta:\theta<\lambda\}$ is not eventually constant. \newline
Then for some $\mu\leq 2^\lambda$ we have Galvin's property on $\mathcal{D}_\lambda$.
\end{corollary}

\par\noindent\emph{Proof}. \newline
Let $\kappa=\bigcup\limits_{\theta<\lambda}2^\theta=2^{<\lambda}$. It follows that $\cf(\kappa)\leq\lambda$ and hence $\kappa<2^\lambda$. By Theorem \ref{inac} (a), Galvin's property holds on $\mathcal{D}_\lambda$ for every collection of $\mu$ club subsets when $\mu\in[\kappa^+,2^\lambda]$.

\hfill \qedref{cc1}

Part $(\aleph)$ of the above question is still open, since Galvin's property might be strictly weaker than the weak diamond.
In \cite{bgh} it is proved that under the assumptions of the above corollary we have $\Psi_\lambda$ ($\Psi_\lambda$ is the very weak diamond principle which says that the set of guesses is unbounded in $\lambda$), and $\Psi_\lambda$ is strictly stronger than Galvin's property.
However, we comment that the above corollary limits the possibility of applying the forcing in \cite{MR830084} to inaccessible cardinals. In fact, we must admit the situation of a weakly inaccessible $\lambda$ so that there are $\theta_0<\lambda<\mu$ for which $2^\theta=\mu$ for every $\theta\in[\theta_0,\lambda]$. We may, therefore, ask:

\begin{question}
\label{q2} Assume $\lambda$ is weakly inaccessible. Is it possible that Galvin's property fails for every $\mu\leq 2^\lambda$?
\end{question}

\newpage

\bibliographystyle{amsplain}
\bibliography{arlist}

\end{document}